\documentclass[12pt]{article}
\input epsf      
\usepackage{graphicx}

\usepackage{amssymb}
\usepackage{amsfonts}
\usepackage{amsthm}
\usepackage{amsfonts}

\newcommand{\disc}{\mathbb{D}}

\newcommand{\boundary}{\partial}
\newcommand{\ra}{\rightarrow}

\newcommand{\ba}[1]{\begin{array}{#1}}
\newcommand{\ea}{\end{array}}

\newcommand{\be}{\begin{equation}}
\newcommand{\ee}{\end{equation}}
\newcommand{\bea}{\begin{eqnarray}}
\newcommand{\eea}{\end{eqnarray}}
\newcommand{\beann}{\begin{eqnarray*}}
\newcommand{\eeann}{\end{eqnarray*}}

\newcommand{\integers}{{\mathbb Z}}

\newcommand{\complex}{{\mathbb C}}

\newcommand{\prob}{{\bf P}}
\newcommand{\cord}{{z}}
\newcommand{\har}{h}

\newcommand{\saws}{{\mathcal W}}
\newcommand{\walks}{{\mathcal S}}

\newcommand {\cent} {{\bf c}}

\def\reff#1{(\ref{#1})}
\begin{document}

\title{
Lattice effects in the scaling limit \\
of the two-dimensional self-avoiding walk
}

\author {
Tom Kennedy \footnote{ University of Arizona; research supported  
by NSF grant DMS-0758649.}
\and
Gregory F. Lawler\footnote {University of Chicago; 
research supported by NSF grant DMS-0907143.}
}

\date{}

\maketitle

\begin{abstract}
We consider the two-dimensional self-avoiding walk (SAW) in a simply
connected domain that contains the origin. The SAW starts at the origin
and ends somewhere on the boundary. 
The distribution of the endpoint along the boundary 
is expected to differ from the SLE partition function prediction 
for this distribution  because of lattice
effects that persist in the scaling limit.
We give a precise conjecture for how to compute this lattice effect correction
and support our conjecture with simulations. We also give a precise 
conjecture for the lattice corrections that persist in the scaling 
limit of the $\lambda$-SAW walk.  
\end{abstract}

\newpage

\section{Introduction}

Let $D$ be a bounded, simply connected domain in the complex plane that 
contains $0$. We are interested in the self-avoiding walk (SAW) in $D$ 
starting at the origin and ending on the boundary of $D$. 
It is defined as follows. We introduce a lattice with spacing 
$\delta>0$, e.g., $\delta \integers^2$.
A self-avoiding walk is a nearest neighbor walk on the lattice with the 
property that it does not visit any site more than once.  
To be precise, let $\saws(D,\delta)$ denote the set of functions
of the form 
 $\omega :\{0,1,2,\cdots,n \} \ra \delta \integers^2$ where
 $n$ is a positive integer;
 $||\omega(i)-\omega(i-1)||=\delta$ 
for $i=1,2,\cdots,n$;  $\omega(i) \neq \omega(j)$ for $0 \le i < j \le n$;
$\omega(0) = 0$; $\omega(j) \in D, j < n$; $\omega(n) \not \in D$.
The integer $n$ is the number of steps in the SAW, and from now on we 
will denote it by $|\omega|$.

Since $\saws(D,\delta)$ is finite, 
 we can define
a probability measure on $\saws(D,\delta)$
by taking the probability of $\omega$ to 
be proportional to $\beta^{|\omega|}$ where $\beta>0$ is a parameter. So 
\bea
\prob(\omega) = \prob_{D,\delta}(\omega) = {\beta^{|\omega|} \over Z(D)}
\eea
where the partition function $Z(D) = Z_\delta(D)$ 
is defined by the requirement 
that this be a probability measure.
One can consider this model for all $\beta>0$, but it is most 
interesting for one particular value that makes the model critical, 
$\beta=1/\mu$, where $\mu$ is the connective constant which we define
next.

Let $c_N$ be the number of SAW's in the lattice with $N$ steps that 
start at $0$. 
(They are not constrained to lie in $D$.) It is known 
that this number grows exponentially with $N$ in the sense that the 
following limit exists \cite{madras_slade}.
\bea
\mu = \lim_{N \ra \infty} c_N^{1/N}.
\eea
The connective constant $\mu$   depends on the lattice.  
Nienhuis \cite{nienhuis} predicted that for the hexagonal lattice
$\mu=\sqrt{2+\sqrt{2}}$, and this was recently proven by 
Duminil-Copin and Smirnov \cite{dc_smirnov}. 
For the square and triangular lattices there are only numerical estimates of 
the value of $\mu$. For the remainder of this paper we will 
take $\beta=1/\mu$ to make the model critical. 

We also consider the analogy of the above definition with the ordinary 
random walk.
 The natural way to describe the random walk in $D$ starting
at $0$ and ending on the boundary of $D$ is to start a random walk 
at $0$ and run it until it hits the boundary.  Let $\walks(D,\delta)$
be the set of such walks.   The probability of a 
particular such random walk $\omega$  is  $\cord^{-|\omega|}$ where $\cord$
is the coordination number of the lattice, e.g., $\cord=4$ for the square
lattice.  We can consider this as a random walk starting at the origin
stopped when it leaves $D$, or just as a measure on $\walks(D,\delta)$
that assigns measure $z^{-n}$ to each walk of length $n$. 
In the SAW the connectivity constant $\mu$ plays the role of the 
coordination number $z$ for the random walk.

In both the SAW and the ordinary random walk we are interested 
in the scaling limit in which the lattice spacing $\delta \ra 0$.
Let us first discuss the case of ordinary random walk which
is well understood.

For the ordinary random walk the scaling limit is Brownian motion 
starting at $0$ and stopped when it hits the boundary of $D$. 
The distribution of the endpoint of the Brownian motion on the 
boundary is harmonic measure.  The lattice effects associated
to the definition of the first boundary point of the lattice
walk disappear in the scaling limit.  The key fact is that if
a random walk or a Brownian motion gets very close to the boundary,
then it will hit it soon.  Therefore, if we couple a random
walk and a Brownian motion on the same probability space so the
paths are close, then the first time that the random walk hits
the boundary will be close to the first time that the
Brownian motion hits the boundary.  See \cite[Section 7.7]{LLimic}
for a precise statement.  The argument there works for any
simply connected domain even with nonsmooth boundaries, and extends
to finitely connected domains as well.  One does need to
assume that  the boundary is sufficiently large so that when
the Brownian motion or random walk gets close, then it is very
likely to hit it soon.

If  the boundary of our domain $D$ is a piecewise smooth curve, 
then harmonic measure is absolutely continuous with respect to arc
length along the boundary \cite{riesz}.
We let $\har_D(z)$ denote its density 
with respect to arc length; this is often called the Poisson
kernel (starting at $0$). If $f$ is a conformal map on $D$ that 
fixes the origin and such that the boundary of $f(D)$ is also 
piecewise smooth, 
then the conformal invariance of Brownian motion \cite{levy}
implies that the density for harmonic measure on the boundary 
of $f(D)$ is related to the density on the boundary of $D$ by
\begin{equation}  \label{jul2.1}
\har_D(z) = |f^\prime(z)| \har_{f(D)}(f(z)).
\end{equation}
If $D$ is simply connected
and we take $g_D$ to be a conformal map of $D$ onto the unit disc which
fixes $0$, then by symmetry $\har_{f(D)}(f(z))$ is just $1/2 \pi$, and so 
$\har_D(z) = |g_D^\prime(z)| / 2 \pi$.   We emphasize that
conformal invariance of harmonic measure implies that the exponent
of $|f^\prime(z)|$  in (\ref{jul2.1}) equals one.

We now
consider the SAW in $D$ from the origin to the boundary of $D$.
Here we will assume that $\partial D$ is piecewise smooth. 
In this case we have the following conjecture.
\begin{itemize}
\item  As $\delta \rightarrow 0$, the measures $\prob_{D,\delta}$
converge to a probability measure $\prob_D$ on simple paths
from the origin to $\partial D$.
\item  The measure $\prob_D$ can be written as
\bea
          \int_{\partial D} \tilde{\rho}_D(z) \, \prob_D(0,z) \,
 |dz| , 
\label{pdref}
\eea
 where $\tilde{\rho}_D(z)$ is the density of a probability measure
 on $\partial D$ and $\prob_D(0,z)$ is the probability measure associated
 to radial $SLE_{8/3}$ from $0$ to $z$ in $D$.
 \item  There exists a periodic function, $l(\theta)$,
such that 
 \bea
\tilde{\rho}_D(z) =  l(\theta(z,D)) \, \rho_D(z)
\label{our_conj}
\eea
where $\theta(z,D)$ is the angle of the tangent to $\partial D$
at $z$ and $\rho_D(z)$ is a multiple of the $SLE_{8/3}$ partition
function. The function $l(\theta)$ and its period depend on the lattice. 
For example, on the square lattice the period is $\pi/2$. 
\item  The density $\rho_D(z)$ transforms under
conformal maps by 
\bea
\rho_D(z) = c |f^\prime(z)|^{5/8} \rho_{f(D)}(f(z))
\label{conf_covar}
\eea
The constant $c$ is determined by the constraint that this be a probability
density. (It depends on $D$ and $f$.) In particular, if
$D$ is simply connected, and $g_D: D \rightarrow \disc$
with $g_D(0) = 0$, then 
\bea
\rho_D(z) = c |g_D^\prime(z)|^{5/8}
\label{conf_covar_g}
\eea
($\disc$ denotes the unit disc centered at the origin.)
\end{itemize}

In the case of $\integers^2$ 
where the boundary of $D$ is composed of horizontal
and vertical line segments, this conjecture was made
by Lawler, Schramm and Werner \cite{lsw_saw}.  Simulations
on an infinite horizontal strip \cite{dyhr_etal} give strong
support to the conjecture.   For such
boundaries, $l(\theta(z,\theta))$ is constant.    The conjecture
was reiterated in \cite{lawler_utah} where it was also 
conjectured for other domains ``...after taking care of the local
lattice effects.''  Our conjecture makes precise
the nature of the lattice effects in terms
of the lattice correction $l(\theta(z,\theta))$.
 The function $l(\theta)$ will depend on 
just how we define ``ending on the boundary of $D$'' and on the lattice type.
In general, we use $\rho$ to denote densities that do not include the lattice
effects that persist in the scaling limit, e.g., \reff{conf_covar}, 
and we use $\tilde{\rho}$ to denote densities that do include the 
lattice effects, e.g., \reff{our_conj}.

Unfortunately, the   Monte Carlo algorithms for simulating the above ensemble 
are local algorithms and so are not very efficient \cite{madras_slade}.
We have not attempted to test the
conjecture for this ensemble by simulation. Instead we introduce
another ensemble that we can study with the pivot algorithm, 
a fast global Monte Carlo algorithm.   Instead of stopping a walk at a
boundary point, one chooses an infinite length walk conditioned on
the event that it crosses the boundary only once.
We will refer to this ensemble as the ``cut-curve'' ensemble since
the boundary of the domain cuts the SAW into two SAW's, one contained 
in $D$ and one in the complement of $D$.    The above ensemble and the cut-curve
ensemble are the same for the infinite strip
studied in \cite{dyhr_etal}, but for most domains this is not the case.

The scaling limits of
 the SAW and the loop-erased walk, which is obtained by erasing
 loops from the ordinary random walk, are two cases of
 the Schramm-Loewner evolution (SLE).  The discrete models can be
 considered as special cases of the $\lambda$-SAW.  We review
 the $\lambda$-SAW in the next section, and we extend our
 conjecture to this case.  (There is no precise conjecture in the
 literature on the nature of the lattice correction, and we think it
 is worthwhile to write it down.)  After that we consider the
 conjectured scaling limits of the two ensembles --- walks stopped
 upon reaching the boundary and walks conditioned to hit the
 boundary only once. In only the $\kappa = 8/3$ (SAW) case do we
 expect an equivalence of these ensembles, and therefore the tests
 we do here would not work for other values of $\kappa$.  

In section two we return to the SAW and give explicit conjectures 
for the lattice correction function $l(\theta)$ for the two
lattice ensembles.  In section three, we discuss simulations for the
cut-curve ensemble including numerical calculation of the lattice correction
function. In the final section we summarize our results 
and discuss the lattice correction function $l(\theta)$ 
for other interpretations of the SAW ending on the boundary of the domain.

\subsection{$\lambda$-SAW}

The conjectures  for the SAW  and the results for the ordinary
walk (considered in terms of the loop-erasure of the paths)
are  particular cases of
conjectures for 
a model called
the $\lambda$-SAW introduced by Kozdron and Lawler \cite{LKozdron}.
Even though we are only testing the SAW conjecture, 
 we will give the conjectures
for the general model.
There are two versions,  chordal (boundary-to-boundary)
and radial (boundary-to-interior);
we will restrict our discussion here to the radial case.
(There is also an interior-to-interior case, but then there is no
boundary lattice correction so it is not relevant for this paper.)
As above, we assume that $D$ is a domain with piecewise smooth
boundary.  For convenience we  use $\integers^2$ for our lattice, but the
definition can be extended to other lattices.  The parameter
$\lambda$ can be considered a free parameter, but we will set
  $\lambda = -\cent/2$ where 
  $\cent \leq 1$ denotes   central charge.  (We will not define
  central charge in this paper and can just take it as a parameter.)
  The $\lambda$-SAW is a model that is conjectured to have
  a scaling limit of $SLE_\kappa$ where 
  \[  \kappa =   \frac{13 - \cent - \sqrt{(13-\cent)^2 - 144}}
                {3} \in (0,4] .\]

The (rooted) random walk loop measure is the measure on ordinary
random walk loops that assigns weight $(1/4)^n \, n^{-1}$ to
each loop of $n >0$ steps.   A loop is a path which begins and ends
at the same point.  The measure can also be considered as 
a measure on unrooted loops by forgetting the root.

For each lattice spacing $\delta$ and each SAW $\omega$ as above,
we let $m_D(\omega)$ denote the total measure of the set of loops
that lie entirely in $D$ and share at least one site with $\omega$.
If $\beta >0$, the $\lambda$-SAW gives each SAW $\omega$
as above weight
\[             q(\omega)
=  q_{\delta,D,\cent,\beta}(\omega)  =  \beta^{|\omega|}\,
  \exp \left\{- \frac {\cent}{2}
   \, m_D(\omega) \right\} . \]
The partition function is
\[    Z(D;\delta) = Z(D;\delta,\cent,\beta) = \sum {q(\omega)} , \]
where the sum is over all SAWs on the lattice $\delta \integers^2$ that
start at the origin and end at the boundary.  (As noted before, there
are several definitions for ``ending at the boundary''.  In the discussion
here, we fix one such definition and the lattice-dependent quantities depend
on the choice.)  It is conjectured that for each $\cent$, there is
a critical value $\beta = \beta_\cent$ such that $Z(D;\delta)$ follows
a power law in $\delta$ as $\delta \rightarrow 0$.  We assume this
conjecture and fix $\beta$ at the critical value.

If $\cent =0$, then this is the usual SAW model.  If $\cent = -2$,
then this is the loop-erased random walk (LERW) which is obtained
by taking the ordinary random walk as above and erasing loops
chronologically from the paths.  The partition function $Z(D;\delta)$
for the loop-erased walk is exactly the same as that for the
usual random walk.  (See \cite[Chapter 9]{LLimic} for a discussion
of this.)  We state our conjectures in terms of the boundary and interior
scaling exponent for $SLE_\kappa$:
\[   
    b = \frac{6-\kappa}{2\kappa} \in [ 1/4,\infty)
    , \;\;\;\; \tilde b = \frac{b\, (\kappa-2)}{4} .
\]

\begin{itemize}

\item  There is a lattice correction function $l(\theta)$ that
is continuous, strictly positive and periodic.

\item  If $\omega$ is from $0$ to $\partial D$, let $l(\omega)$
denote $l(\theta(z,D))$ where $z$ is the first point on $\partial D$
hit by a bond of $\omega$ and $\theta(z,D)$ is the angle of the
tangent to $\partial D$ at $z$.  Define the lattice-corrected
weight by
\[             \hat q(\omega) = \frac{ q(\omega)}{ l(\omega)}.\]

\item As $\delta \rightarrow 0$, the measure $\delta^{1-b-\tilde b}
 \, \hat q$ on paths converges to a nontrivial finite
 measure $\nu_D$ on simple paths from $0$ to $\partial D$.  It
 can be written as
 \[             \nu_D = \int_{\partial D} \rho_D(z) \, \nu_D^\#(0,z)
  \, |dz| , \]
  where $\rho_D$ is a positive function and $\nu_D^\#(0,z)$ is
  a probability measure on simple paths starting at $0$ and leaving
  $D$ at $z$.  
  
  \item  If $g$ is a conformal transformation with $g(0) = 0$
   that is smooth
  on $\partial D$, then 
  \begin{equation}  \label{jul3.1}
             \rho_D(z) = |g'(z)|^b \, |g'(0)|^{\tilde b}
   \, \rho_{g(D)}(z).
   \end{equation}
   
   \item  The probability measures $\nu_D^\#(0,z)$, considered as
   measures on curves modulo reparametrization, are conformally invariant.
   More specifically, if $\gamma:[0,t_0) \rightarrow D$ is a curve
   with $\gamma(t_0-) \in \partial D$,  let
   \[             \sigma(t) = \int_0^t |g'(\gamma(s))|^d \, ds, \;\;\;\;
     d = 1 + \frac \kappa 8 ,\]
   and define $g \circ \gamma(t)$ by
   \[             g \circ \gamma(\sigma(t)) = g(\gamma(t)), \;\;\;\;
     0 \leq t \leq \sigma(t_0). \]
If $g \circ \mu^\#_D(0,z)$ denotes the induced measure on curves on $g(D)$,
then 
\[             g\circ \nu_D^\#(0,z) = \nu_{g(D)}^\#(0,z) . \]

\item  If $D$ is simply connected, then $\nu_D^\#(0,z)$ is the reversal of 
radial $SLE_\kappa$ from $z$ to $0$ in $D$ with the natural parametrization.
(This should also be true for multiply connected $D$ under the appropriate
definition of $SLE_\kappa$ proposed in \cite{lawler_partfunc}.)

\end{itemize}

These conjectures are a long way from being proved. Indeed, a special
case is the SAW model which is a notoriously difficult problem!  However,
we state them here to see that the precise conjecture requires discussing
a boundary lattice correction; our conjecture is that the correction only
depends on the angle of the boundary.   We make a number of comments.

\begin{itemize}

\item  The density $\rho_D(z)$ is sometimes called the partition function
for radial $SLE_\kappa$.  It is defined up to an arbitrary constant.

\item  The loop-erased walk ($\cent = -2$) is particularly nice because
it is closely related to the ordinary random walk.  The partition function
is the Poisson kernel (density of harmonic measure)
even for non-simply connected domains.  

\item  For other values of $\cent$, if $D$ is simply connected,
the partition function is
the Poisson kernel raised to the $b$th power.  However, this is not
true for multiply connected domains.

\item  The loop terms $m_D(\omega)$ can be considered as having three parts.
The very small loops that occur away from the boundary contribute
a microscopic (lattice dependent) part that affects the
critical value of $\beta$.  
The large loops contribute a macroscopic term that is
seen in the scaling limit; this is the Brownian loop measure as
defined in \cite{LWsoup}.  Finally,
there are the small loops that occur near the boundary.  They give
both a macroscopic effect seen in  the exponent $b$ and a microscopic
  effect in the function $l$.  The boundary effect is measured both
in the number of walks that stay one one side of a line and in the
measure of loops that go on the walks.  We only see the first effect
in the SAW ($\cent = 0$) case.

\item  If $\kappa = 8/3$, then $b=5/8, \tilde b = 5/48$.  There is a
$|g'(0)|^{\tilde b}$ factor in \reff{jul3.1} that does not seem to appear
in \reff{conf_covar}.  However, it is implicitly there in the normalization
to make the measure a probability measure.

\end{itemize}

\subsection{Cut-curve configurations}

We will be testing a cut-curve ensemble for SAW's.  The case
$\cent = 0$, which is what we use in this paper,
 is special for such configurations and agrees with
the bridge decomposition of restriction measures \cite{albert_hugo},
but for the sake of completeness let us discuss
the general case for $\cent \leq 1$.  Suppose $D$ is a bounded,
simply connected domain containing the origin whose boundary
$\partial D$ is a smooth Jordan curve.  Let $D^*$ be the unbounded component of
$\complex \setminus \partial D$.  We consider two measures on
simple curves from $0$ to infinity that intersect $\partial D$ only
once:
\begin{enumerate}
\item  Take ``whole-plane'' $SLE_\kappa$ and condition on the event that
the curve intersects $\partial D$ only once.  This is conditioning on an event
of probability zero, so one must define this in terms of a limit.
\item  Take independent copies of radial $SLE_\kappa$ in $D$ and in $D^*$
and condition them to hit the same point in $\partial D$.  Then concatenate
the paths.
\end{enumerate}
In order to see the difference, let us describe the lattice models that
we expect to converge to these measures.  Since it is hard to talk
about infinite walks, we will choose a point $z^* \in D^*$ with large
absolute value.  We take the scaling limit for fixed $z^*$ and then
let $z^*$ go to infinity.  If $\delta$ is a scaling factor, we abuse notation
slightly and write $z^*$ for the point in $\delta \integers^2$ closest to $z^*$. 
 We consider two measures on paths.  
 For each $\delta$, we consider the set $\saws(\delta,D,z^*) $
 of SAW's $\omega$
on $\delta \integers^2$ with $\omega(0) = 0, \omega(|\omega|) = z^*$,
$|\omega(j)| 
\leq |z^*|^2$ for all $j$
such that there exists only one bond 
that intersect $\partial D$.  We write $z_\omega, w_\omega$ for
the vertices in this bond with $z_\omega \in D, w_\omega \in \complex
\setminus D$, and we write $\omega = \omega^D \oplus 
[z_\omega,w_\omega] \oplus \omega^{*}$ 
where
\[ \omega^D = [\omega(0), \ldots,z_\omega], \;\;\;\;
  \omega^{*} = [w_\omega,\ldots,\omega(|\omega|)]. \]
Let $\beta = \beta_\cent$ be the critical value.  Then we consider
the following measures on $\saws(\delta,D,z^*) $:
\begin{enumerate}
\item  Each $\omega$ gets weight 
\[ q_1(\omega) = \beta^{|\omega|}
 \, \exp\left \{- \frac \cent 2 m(\omega;z^*)\right\}, \]
where $m(\omega;z^*)$ denotes the measure of loops staying in 
the ball of radius $|z^*|^2$ that intersect $\omega$.
\item  Each $\omega$ gets weight
\[  q_2(\omega)  =  \beta^{|\omega|}
 \, \exp \left\{- \frac \cent 2 [m_D(\omega^D)+ m_{D^*}(\omega^{*})]\right\} \]
 \end{enumerate} 
 We can state our conjectures as follows.  Let $\theta(\omega)$ be the
 angle of the tangent to $\partial D$ at the point where $\omega$ hits
 $\partial D$. 
 \begin{enumerate}
 \item  There exists a lattice correction function $\hat l_1$ such that we
 can take the 
 scaling limit of the measure
 \[             \hat q_1(\omega) = \frac{q_1(\omega)}{\hat l_1(\omega)} . \]
 If we then take $z^* \rightarrow \infty$, we get whole plane $SLE_\kappa$
 conditioned to hit $\partial D$ only once.  
 \item  There exists a lattice correction function $\hat l_2$ such that the
 scaling limit of the measure
 \[    \hat q_2(\omega) = \frac{q_2(\omega)}{\hat l_2(\omega)}\] exists.
 If we then take $z^*$ to infinity, we get the measure 
 given by
 \[        c \int_{\partial D} \rho_D(z) \, \rho_D^*(z) \, 
    \left[\nu_D^\#(0,z) \oplus  \nu_{D^*}^\#(z,\infty) \right] \, |dz| . \]
   Here $\rho^*(z)$ is the density as in  \reff{jul3.1} for radial $SLE_\kappa$
   in $D^*$ centered at infinity. 
\end{enumerate}

For $\cent = 0$ (restriction measures), the two limits agree and this
is what we use for SAW in this paper.  For other values of $\cent$ we
get different measures.  For example if $\cent = -2$ (loop-erased walk),
the first measure corresponds to loop-erased walk conditioned to hit
$\partial D$ only once and the second measure corresponds to the loop-erasure
of an ordinary walk conditioned to hit $\partial D$ only once.

\section{Lattice effects}

The constraint that the SAW stays in $D$ has both a macroscopic 
and a microscopic effect on the boundary density. 
The macroscopic effect is captured by the 
 conjecture \reff{conf_covar}.
The microscopic effect comes from the behavior near the endpoint of the walk 
on the boundary of $D$. 
Consider a SAW that ends at $z \in \boundary D$ 
and consider the tangent line to the boundary at $z$. 
The constraint that the SAW stays in $D$ implies that near $z$ 
the SAW must stay on one side of this line. 
Loosely speaking, the number of SAW's that end at $z$ and stay on one side
of this line depends on the orientation of the line with respect 
to the lattice. The result is a factor $l(\theta)$ that depends on 
the angle of the tangent line with respect to the lattice, 
and so we obtain our conjecture \reff{our_conj}.

The lattice correction function $l(\theta)$ depends on the type 
of lattice and on how we interpret ``ending on the boundary of $D$.''
We will first discuss the interpretation we introduced at the 
start of the introduction. 
We consider all SAW's $\omega$ with 
$\omega(0)=0$, $\omega(i) \in D$ for $i=0,1,2,\cdots,|\omega|-1$ and 
$\omega(|\omega|) \notin D$. ($|\omega|$ denotes the number of steps
in $\omega$.) So the last bond of the SAW intersects the boundary of $D$, 
and this is the only bond in the SAW that intersects the boundary.  

We can compute the lattice correction function $l(\theta)$ as follows. 
We give the details for the square lattice. 
Other lattices, e.g.,  the triangular or  hexagonal, 
will require some modifications.
Consider a bond that intersects the boundary $\boundary D$. Let $z$
be the endpoint of the bond that is in $D$. 
Let $w$ be the point where the bond intersects the boundary
(typically not a lattice site).  
Consider the tangent line to the boundary at $w$. 
We need to count the number of SAW's of a fixed length $N$ that 
end at $z$ and do not intersect this tangent line. This quantity will 
depend on the angle of the tangent line with respect to the lattice. 
It will also depend on the distance $|w-z|$ and on whether the bond 
is horizontal or vertical. However, these two factors are the only 
dependence on the bond. 

Motivated by the above, we consider the vertical bond between $(0,0)$ 
and $(0,1)$. 
We fix an $l$ in the interval $[0,1]$ and an angle $\theta$.
The parameter $l$ plays the role of the distance $|w-z|$.  
Let $L$ be the line with polar angle $\theta$ passing through the point 
$(0,l)$. We consider SAW's with $N$ steps ending at the origin which
by reversal can be considered as beginning at the origin.
Let $c_N$ be the number of such walks, and let
$a_N(l,\theta)$ be the number of such walks that do not intersect the line.
So $a_N(l,\theta)/c_N$ is the probability that the
$N$ step SAW does not intersect the line. We expect that
there exists a function $p^v(l,\theta)$ such that 
\beann
{a_N(l,\theta) \over c_N} \sim  p^v(l,\theta) \, N^{-\rho}, \quad N \rightarrow \infty
\eeann
with $\rho=25/64$ \cite{lsw_saw}.  The actual exponent is not important here.
The key is that if $l', \theta'$ are two different values, then
\[            \frac{a_N(l,\theta)}{a_N(l',\theta')} \sim
   \frac{p^v(l,\theta)}{p^v(l',\theta')},\quad N \rightarrow \infty.\]
   
We define 
\bea
p^v_1(l,\theta)= \lim_{N \ra \infty} {a_N(l,\theta) \over c_N} N^{\rho}
\eea
(We do not know how to prove this limit exists.)
The superscript $v$ on $p$ indicates that we 
took the bond crossing the boundary to be vertical.
We use the subscript $1$ to distinguish the quantities related 
to the lattice correction function for this particular ensemble from the 
ensemble that we will consider next.
We let $p^h_1$ denote the analogous quantity for a horizontal bond. 
The symmetry of the square lattice implies that
$p_1^h(l,\theta)=p^v_1(l,\theta+\pi/2)$. (Note that $p^h_1(l,\theta)$ and 
$p^v_1(l,\theta)$ have period $\pi$ in $\theta$.)

Now consider what happens as we move along the boundary. The angle $\theta$
of the tangent will vary smoothly. The distance $l$ will not. 
As long as $\tan(\theta)$ is not rational, $l$ will be distributed
uniformly between $0$ and $1$. 
So in the scaling limit, averaging over an infinitesimal section of 
the boundary will be equivalent to averaging $l$ over $[0,1]$. So we define 
\beann
p_1^x(\theta) = \int_0^1 p_1^x(l,\theta) \,dl,  \quad x = h,v
\eeann

The function $p^x_1(\theta)$ captures the microscopic lattice effect 
caused by the constraint that the SAW must stay on one side of 
a line as it approaches the boundary. There is another lattice effect 
that comes from the density of bonds that cross the boundary.
Define $b^v(\theta)$ to be the density of vertical bonds along a line 
with polar angle $\theta$, i.e., the average number of vertical 
bonds that intersect the line per unit length. 
We have $b^v(\theta)=|\cos(\theta)|$. The density of horizontal 
bonds $b^h(\theta)$ is $b^v(\theta-\pi/2)=|\sin(\theta)|$.
The lattice correction function is then 
\bea
l_1(\theta) = b^v(\theta) p^v_1(\theta) + b^h(\theta) p^h_1(\theta)
\label{ldef}
\eea
and the boundary density for this first interpretation of ending on the 
boundary is 
\bea
\tilde{\rho}_{D,1}(z) = \tilde{c_1} |g_D^\prime(z)|^{5/8} l_1(\theta(z,D))
\eea

\medskip

The Monte Carlo algorithms for simulating the above ensemble 
are local algorithms and so are not very efficient \cite{madras_slade}.
We have not attempted to test the
conjecture for this ensemble by simulation. 
Instead we use the cut-curve ensemble introduced earlier.  
We consider the infinite length SAW in the full plane starting at $0$. 
We condition on the event that the SAW crosses the boundary of $D$
exactly once. Consider a bond crossing the boundary, and let $z$ 
be the endpoint inside $D$ and $w$ the endpoint outside of $D$. 
If we condition on the event that the SAW contains this particular bond
and this is the only bond in the SAW that 
crosses the boundary, then we have a SAW
in $D$ from $0$ to $z$ and a SAW in the exterior of $D$ from $w$ 
to $\infty$. So the   conjecture for the 
boundary density will be a product of two functions. The 
interior SAW from $0$ to $z$ gives a factor of 
$\rho_D(z)$ with $\rho_D(z)$ given by \reff{conf_covar}.
The exterior SAW from $w$ to $\infty$ gives a factor of 
$\rho^*_D(z)$ given by
\beann
\rho^*_D(z)= c |h^\prime_D(z)|^{5/8}
\eeann
where $h_D(z)$ is the conformal map of $D^*$ onto the unit disc $\disc$
with $h_D(\infty)=0$. 
Our conjecture for the density of the point where the SAW crosses the 
curve is 
\bea
\tilde{\rho}_{D,2}(z) = \rho_D(z) \, \rho^*_D(z) \, l_2(\theta(z,D))
\label{conj_cut_curve},
\eea
where $l_2(\theta)$ is the lattice correction function 
for the cut-curve ensemble. (The subscript $2$ indicates that the 
quantities are for the cut-curve ensemble.)

We can compute $l_2(\theta)$ as follows.
Again, we restrict our attention to the square lattice.
For $l \in [0,1]$ and an angle $\theta$, let $L$ be the line with 
polar angle $\theta$ passing through $(0,l)$. 
We consider SAW's with $N=2n+1$ steps such that the middle bond
is the bond between $(0,0)$ and $(0,1)$. Let $d_N$ be the number of 
such walks, and let $b_N(l,\theta)$ be the number of such walks that 
only intersect the line once. (Of course it must be the middle 
bond that has the intersection).  
(By translating our SAW's so that they start at $0$, we see that $d_N=c_N/2$.)
Note that $b_N(l,\theta)/d_N$ is the probability the $N$ step SAW 
has no intersection with the line other than the middle bond.
We can think of generating this $N$ step SAW by first generating
two $n$ step SAW's which are independent and attach to $(0,0)$ and $(0,1)$
and then keeping them only if they mutually self avoid. 
The probability they 
are mutually self-avoiding goes as $N^{1-\gamma}$ with $\gamma=43/32$. 
The probability that both of the $n$ step SAW's do not intersect the line
goes as $N^{-2 \rho}$ with $\rho=25/64$ \cite{lsw_saw}.
Note that if both of them do not intersect the line, then they are 
mutually self-avoiding. Hence we expect that there exists 
a function $p^v_2(l,\theta)$ such that 
\beann
{b_N(l,\theta) \over d_N} \sim 
p^v_2(l,\theta) \, N^{-2 \rho + \gamma-1}
\eeann

We define
\beann
p^v_2(l,\theta)= \lim_{N \ra \infty} {b_N(l,\theta) \over d_N} 
N^{2 \rho - \gamma +1}. 
\eeann
Then we define 
\beann
p^v_2(\theta) = \int_0^1 p^v_2(l,\theta) dl
\eeann
As in our derivation of the lattice correction function 
for the first ensemble, the integral over $l$ comes from 
averaging over bonds crossing a small section of the boundary.
We define $p^h_2(l,\theta)$ and $p^h_2(\theta)$ in the analogous way.
As before we let $b^h(\theta)$ and $b^v(\theta)$ denote 
the densities of horizontal and vertical bonds along a line 
with polar angle $\theta$.
The lattice correction function is then 
\bea
l_2(\theta) = b^v(\theta) p^v_2(\theta) + b^h(\theta) p^h_2(\theta)
\label{ldef_cut}
\eea

In the above discussion we have considered SAW's starting at an
interior point in the domain. The same discussion applies to an ensemble
of SAW's that start at a prescribed boundary point of $D$. 
Th  conjecture for the boundary density 
again transforms according to 
\reff{conf_covar}.
With the starting point on the boundary this density is not normalizable. 
We must restrict the endpoint of the SAW to a subset of the boundary 
that is bounded away from the starting point to get a normalizable density. 
A useful reference domain in this case is the upper half plane with the 
starting point at $0$. The unnormalized density for the harmonic measure is 
$x^{-2}$ and so for the SAW it is $x^{-5/4}$. 

\section{Simulations} 

In this section we study the cut-curve ensemble by Monte Carlo 
simulations. There are two types of simulations.
We compute the lattice correction function $l_2(\theta)$ by simulation, 
and we do simulations of the SAW in two different
geometries to test the conjecture \reff{conj_cut_curve}.
We first discuss the computation of $l_2(\theta)$. 

Recall that for odd $N$, $d_N$ is the number of SAW's with $N$ 
steps such that the middle bond is the bond between $(0,0)$ and $(0,1)$. 
For $l \in [0,1]$ and an angle $\theta$, 
$b_N(l,\theta)$ be the number of such SAW's whose only intersection 
with the line through the point $(0,l)$ at angle $\theta$ is in this
middle bond.
The ratio $b_N(l,\theta)/d_N$ is a probability and so may be computed
as follows. We use the pivot algorithm to generate SAW's with $N$ 
steps that start at the origin and such that the middle bond is vertical.
We then pick a point on this bond uniformly at random and take the 
line with angle $\theta$ to go through this point. 
We test if the only intersection of the SAW with the 
line is through the middle bond. 
We find the fraction of the samples that satisfy this condition 
and multiply it by $N^{2 \rho - \gamma +1}$. The result is an 
estimate of $p^v_2(\theta)$. Note that we have included the integral
over $l$ from $0$ to $1$ in the simulation by randomly choosing
$l$ uniformly from $[0,1]$ for each sample. 

\begin{figure}[tbh]
\includegraphics{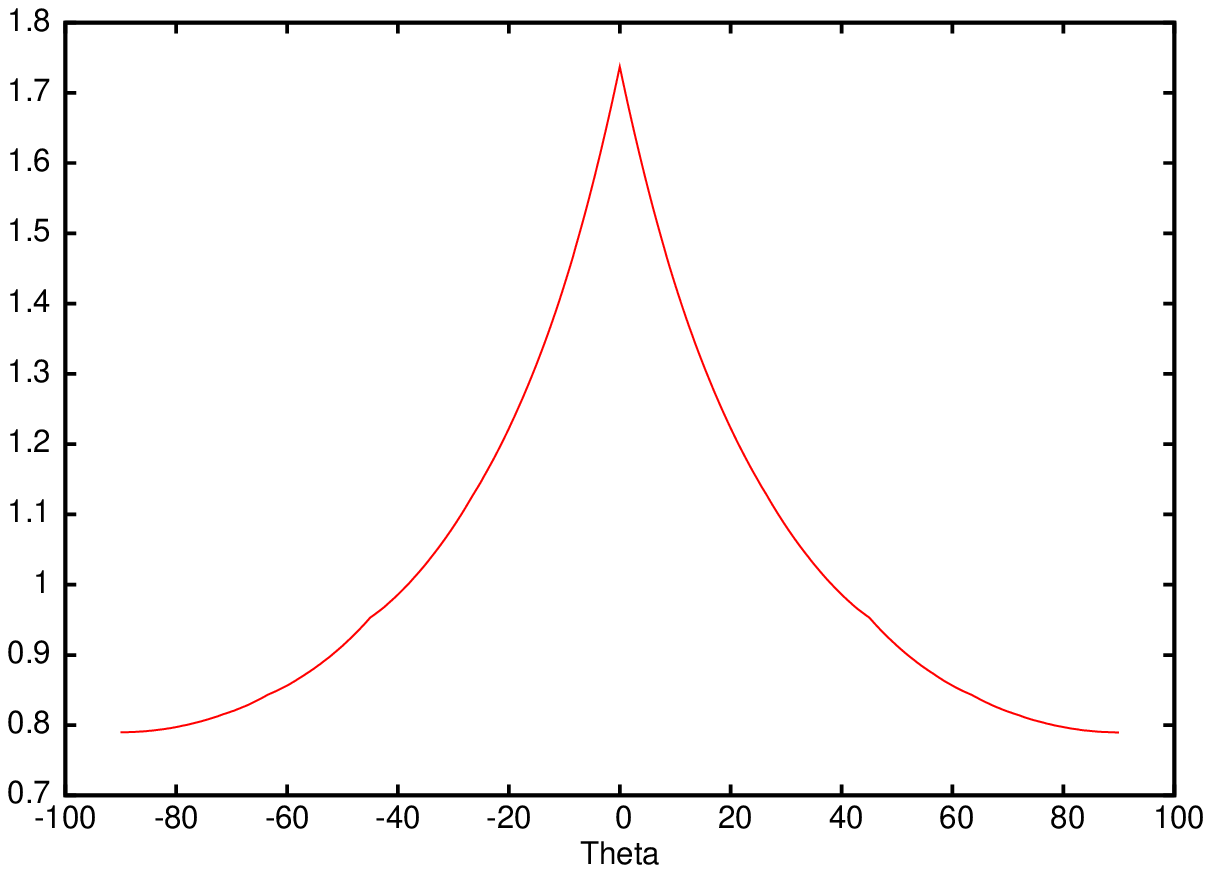}
\caption{\leftskip=25 pt \rightskip= 25 pt 
The function $p^v_2(\theta)$. 
}
\label{pfig}
\end{figure}

We did this simulation for values of $N$ ranging from $101$ to $5001$. 
For the smaller values of $N$, one can clearly see finite $N$ effects.
As is always the case with simulations of the SAW, the time required grows
with $N$. However, in this simulation this is exacerbated by 
the fact that the probability the SAW only intersects the line once goes 
to $0$, and we must multiply the probability we compute in the simulation 
by $N^{2 \rho - \gamma +1}$. 
Even with $2$ billion samples, our simulations for $N=2001$ and 
$5001$ have significant statistical errors. The results for $N=501$ and 
$N=1001$ differ by at most $0.05\%$, and 
we use $N=1001$ in this paper. We generated 
$5$ billion samples for this case which took about 67 cpu-days on 
a rather old cluster with 2.4 GHz cpu's. 
Figure \ref{pfig} shows the function $p^v_2(\theta)$. 
In all our figures we give the angle $\theta$ in degrees. 

We now turn to the simulations to test conjecture \reff{conj_cut_curve}.
We use the pivot algorithm to generate walks in the full plane 
with a constant number of steps $N$. We take $N=1,000,000$.  
The lattice spacing and $N$ are such 
that the size of the SAW is large compared to the domain $D$ so that the 
SAW is effectively infinite. We condition on the event that the SAW 
intersects the boundary of $D$ exactly once. The probability of this 
event goes so zero in the scaling limit, so we must generate very large
numbers of SAW's to get good statistics. 
We use Clisby's fast implementation of the pivot algorithm using 
binary trees \cite{clisby}.

In our first test we take the domain for the cut-curve ensemble 
to be a disc centered at $0$ where the SAW starts. 
We take the lattice spacing to be $N^{-\nu}$ and take the radius of the 
disc to be $R=0.2$. (With $R=0.3$ the effect of the finite length of 
the SAW begins to be noticeable. At $R=0.4$ it is quite noticeable.)
In the simulation we sample the Markov chain every $1000$ iterations
and generate a total of approximately $47$ million samples. 
Just over $10 \%$ of these samples satisfy the condition that the 
SAW only intersects the boundary of the circle once, and  
we have approximately $4.9$ million samples of the boundary density. 

\begin{figure}[tbh]
\includegraphics{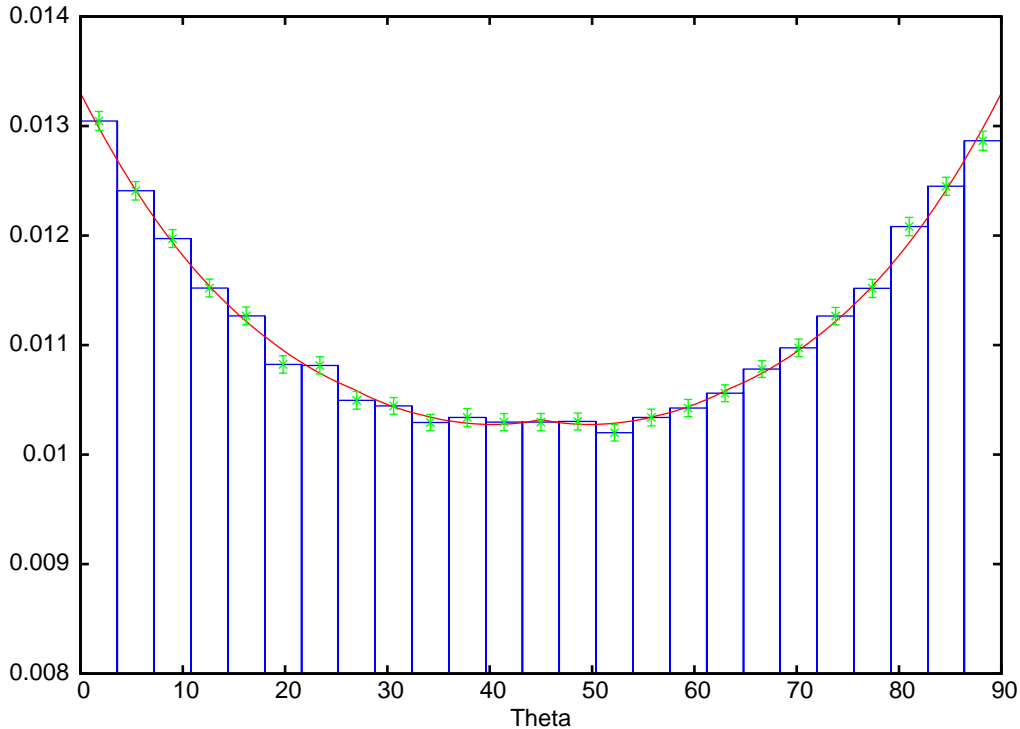}
\caption{\leftskip=25 pt \rightskip= 25 pt 
The function $l_2(\theta)$ (red curve) and the density for the 
full plane SAW conditioned to hit a circle exactly once (blue histogram).  
Note that the range of the vertical axis does
not start at $0$. 
\label{lat_func}
}
\end{figure}

In this geometry both $\rho_D(z)$ and $\rho^*_D(z)$ are constant, 
so the prediction for the boundary density without lattice effects is 
just the uniform density. 
The angle $\theta(z,D)$ of the tangent line at $z$ is equal to the 
polar angle $\theta$ of $z$ mod $90$ degrees. 
So if we think of the boundary density as a function of the polar angle
$\theta$, then our conjecture \reff{conj_cut_curve}
is that the boundary density is proportional to 
$l_2(\theta)$. Figure \ref{lat_func} shows the function 
$l_2(\theta)$ and the boundary density 
we find in the simulation of the cut-curve ensemble. 
Both functions have a period of $90$
degrees. We plot the boundary density as a function of $\theta$ mod $90$. 
Both curves are normalized so that the total area under the curves is one.

Figure \ref{lat_func} compares densities. 
Actually, the function plotted for the simulation of the cut-curve 
ensemble is a histogram. So the points plotted correspond to the average
value of the density over a small interval. 
The simulations do not compute densities directly.
Finding the density requires taking a numerical derivative, i.e., 
computing a histogram. 
We can avoid this extra source of numerical error by working with  
cumulative distribution functions (cdf's) rather than densities.  

\begin{figure}[tbh]
\includegraphics{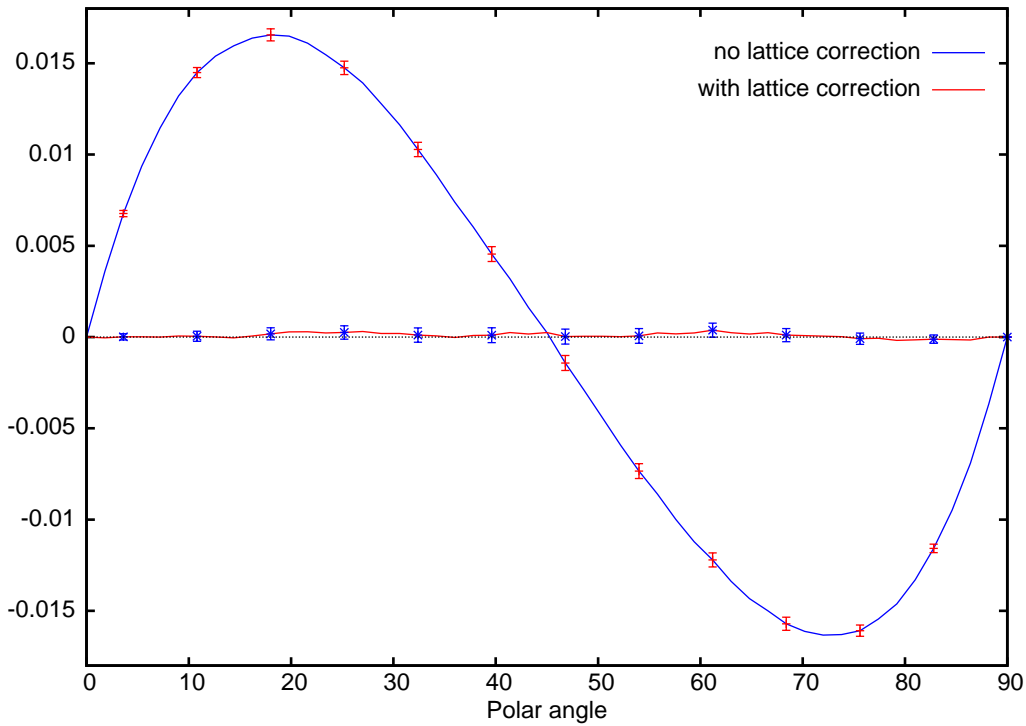}
\caption{\leftskip=25 pt \rightskip= 25 pt 
The simulation of the cdf for the boundary density for the cut-curve 
ensemble for the first geometry minus the theoretical prediction.
The larger curve does not include the lattice correction function; the 
smaller curve does.
}
\label{circle_cdf_per}
\end{figure}

In figure \ref{circle_cdf_per} we study the cdf for the cut-curve ensemble 
using the unit disc centered at $0$. We plot two curves. One is the 
cdf we find in the simulation of the cut-curve ensemble minus the 
cdf for the uniform density, i.e., the density given by 
\reff{conf_covar_g}.
(In this figure we have again taken advantage of the periodicity 
of the underlying density functions.)
This difference is small with the maximum being slightly less than
$2 \%$, but it is clearly not zero. 
In the second curve we show the cdf for the simulation of the 
cut-curve ensemble minus the cdf corresponding to the density 
given by \reff{conj_cut_curve}, i.e., corresponding to $l_2(\theta)$. 
The difference is on the order of $0.02\%$.
The error bars in the figure are two standard deviations for the 
statistical errors, i.e., the error that comes from 
not running the Monte Carlo simulation forever.
There are also errors in the simulation from two other sources - 
the finite length of the SAW and the nonzero lattice spacing. 
We have studied the error from the finite length of the SAW by 
simulating the ensemble with several values of the radius of 
the disc. With $R=0.2$ 
we believe that the error from the finite length of the SAW is much 
smaller than the statistical errors. The nonzero lattice spacing 
means that all our random variables are at a small scale discrete 
random variables. This is reflected in the slightly jagged nature 
of the second curve. 
The error from the nonzero lattice spacing appears to be comparable 
in size to the statistical error. 
Thus the difference between the cdf from the simulation 
and the cdf given by \reff{conj_cut_curve} is 
zero within the errors in our simulation. 
Figure \ref{circle_cdf_per} gives evidence that there are indeed lattice
effects that must be taken into account in the boundary density and 
that our conjecture correctly accounts for these lattice effects.

For our second test of conjecture \reff{conj_cut_curve},
we consider the SAW in the upper half plane, starting at the origin. 
We take the cut-curve to be a semi-circle centered at the origin. 
Again, we take the lattice spacing to be $N^{-\nu}$ and the radius of the 
semi-circle disc to be $R=0.2$. 
We sample the Markov chain every $1000$ iterations
and generate a total of approximately $27$ million samples. 
Approximately $13 \%$ of these samples satisfy the condition that the 
SAW only intersects the boundary of the circle once, 
and we have approximately $3.6$ million samples of the boundary density 
for this geometry. 

The ensemble consists of all SAW's in the upper
half plane which start at $0$ and only cross the semicircle once. 
In this geometry the arc length along the semicircle equals the polar
angle $\theta$. So we will express densities as functions of $\theta$. 
A simple computation using the conformal map $z+1/z$ shows that the 
interior density $\rho_D(\theta)$ is $[\sin(\theta)]^{5/8}$. 
The exterior density $\rho_D^*(\theta)$ 
is exactly the same. (This is just a consequence of the symmetry of 
our geometry under the inversion $z \rightarrow -1/z$.) 
So our conjecture for the density along the cut-curve is 
proportional to 
\bea
[\sin(\theta)]^{5/4} \, \, l(\theta)
\label{sin_density}
\eea

\begin{figure}[tbh]
\includegraphics{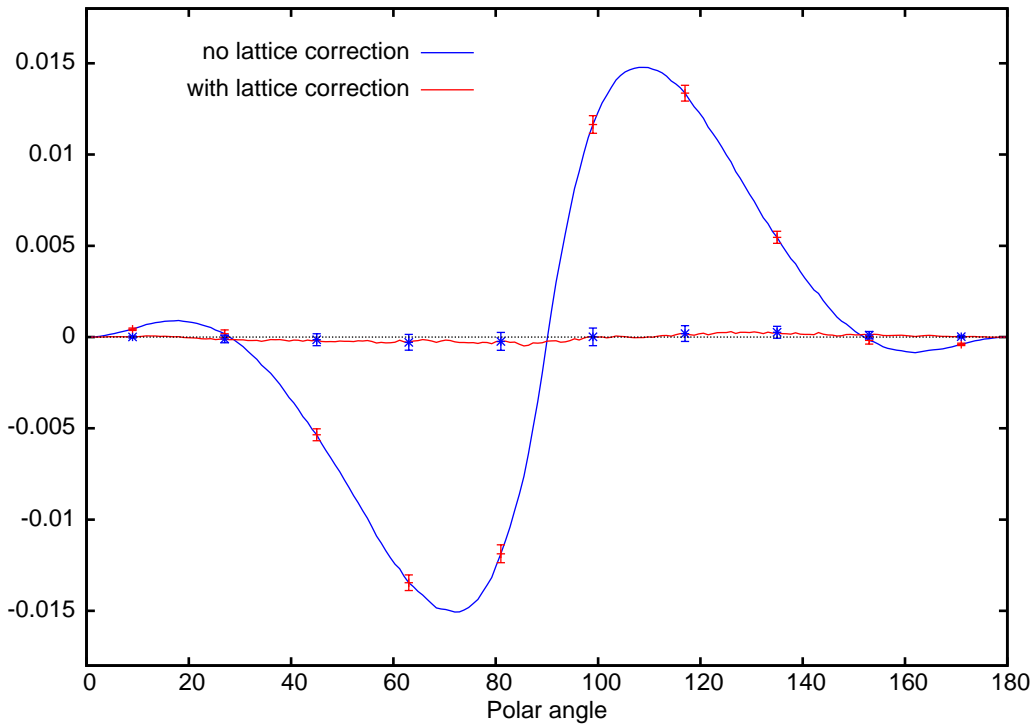}
\caption{\leftskip=25 pt \rightskip= 25 pt 
The simulation of the cdf for the boundary density for the cut-curve 
ensemble for the second geometry minus the theoretical prediction.
The larger curve does not include the lattice correction function; the 
smaller curve does.
\label{semi_dif}
}
\end{figure}

The comparison of our simulation of the cut-curve ensemble cdf for 
the SAW and our conjecture is shown in figure \ref{semi_dif}.
Again we plot two curves. For one curve we find the cdf corresponding to 
just the density function $[\sin(\theta)]^{5/4}$. This would be the 
conjectured cdf if there were no lattice effects. We 
plot the cdf for the simulation of the cut-curve ensemble minus the 
cdf corresponding to just $[\sin(\theta)]^{5/4}$.
This difference is small (the max is on the order of $1.5 \%$), but is 
clearly not zero. For the other curve we compute the cdf corresponding 
to our conjecture \reff{sin_density} with the lattice effect and 
subtract this function from the simulation of the cdf for the cut-curve
ensemble. The difference is on the order of $0.05 \%$ 
which is zero within the errors in our simulation. 
This figure gives further evidence that there are indeed lattice
effects that must be taken into account in the boundary density and 
that our conjecture correctly accounts for these lattice effects.

\section{Conclusions and future work}

We have considered the ensemble of SAW's 
in a simply connected domain containing the origin which start at the 
origin and end on the boundary.  
It has been noted before that the   conjecture 
for this boundary density will have lattice effects that persist 
in the scaling limit. We have conjectured that this lattice effect 
is given by multiplying the   density by a function
$l(\theta(z,D))$ where $\theta(z,D)$ is the angle of the tangent line
to the boundary of $D$ at the point $z \in \boundary D$. The lattice
correction function $l(\theta)$ depends on the lattice and on how we 
interpret ``ending on the boundary of $D$.'' 
We have shown how to compute the lattice correction function $l(\theta)$ 
for two particular interpretations. 
Our focus has been on the distribution of the endpoint of the SAW 
on the boundary, but we should remark that in light of \reff{pdref} 
the lattice effects in this boundary density will produce lattice 
effects in the probability measure $\prob_D$.
We have also extended this conjecture to the $\lambda$-SAW.

As we have noted before, there is no efficient way to simulate the 
natural interpretations of the ensemble of SAW's  in a domain 
which start at an interior point and end on the boundary.
We have circumvented this difficulty by introducing 
the cut-curve ensemble which can be thought of as an ensemble of 
two SAW's, one from the interior point to the boundary  and the other 
from that boundary point to $\infty$. There is another ensemble that
is amenable to efficient simulation which is studied in \cite{tk_dilation}.
Given a domain $D$ containing the origin the ensemble is defined as follows.
We assume the domain has the property that a ray from the origin only 
intersects the boundary of the domain in one point. 
For a SAW $\omega$ we let $\lambda(\omega)>0$ be such that 
the endpoint of $\omega$ is on the boundary of $\lambda(\omega) D$. 
In general the SAW need not be inside the dilated domain $\lambda(\omega) D$. 
Our ensemble consists of all SAW's $\omega$ of any length such that 
$\omega$ is contained in $\lambda(\omega) D$. (One must introduce
cutoffs to make this a finite measure.) In \cite{tk_dilation} we 
show how one can simulate this ensemble using the ensemble of 
SAW's of a fixed length. 

Finally, it is natural to ask if there is an interpretation of 
``ending on the boundary of $D$'' for which there are no lattice
effects in the scaling limit. We speculate that the following ensemble 
has this property. As before, let $\delta$ be the lattice spacing. 
Let $\epsilon>0$ and consider all SAW's that start at the origin, 
stay inside $D$ and end within a distance $\epsilon$ of the boundary 
of $D$. We let $\delta$ go to $0$ first and then let $\epsilon$ go
to zero. We conjecture that $l(\theta)$ is constant for this ensemble. 
Unfortunately, the double limit involved in this ensemble makes it 
difficult to simulate this ensemble.

\end{document}